\newlength{\picturewidth}
\documentclass{article}
\usepackage{amsmath,amsfonts,amssymb,latexsym,graphicx}
\newcommand{\Extra}[1]{}

\allowdisplaybreaks[4]

\newcommand{\Vladimir}{Vladimir}

\newcommand{\DOT}{.}

\newcommand{\st}{\mathrel{:}}

\newcommand{\D}{\,\mathrm{d}}

\newcommand{\true}{^\mathrm{true}}

\newcommand{\bbbp}{\mathbb{P}}		% auxiliary (probability)
\DeclareMathOperator{\Prob}{\bbbp}

\newenvironment{remark*}
  {\trivlist\item[\hskip\labelsep{\bfseries Remark}]\relax}
  {\endtrivlist}

\title{Strong confidence intervals for autoregression}

\author{Vladimir Vovk\\
\texttt{vovk{\rm@}cs.rhul.ac.uk}\\
\texttt{http://vovk.net}}

\begin{document}
\maketitle
\begin{abstract}
  In this short preliminary note
  I apply the methodology of game-theoretic probability
  to calculating non-asymptotic confidence intervals
  for the coefficient of a simple first order scalar autoregressive model.
  The most distinctive feature of the proposed procedure
  is that with high probability it produces confidence intervals
  that always cover the true parameter value
  when applied sequentially.
\end{abstract}

\section{Introduction}
\label{sec:introduction}

Game-theoretic probability
(see, e.g., \cite{shafer/vovk:2001},
with the basic idea going back to Ville \cite{ville:1939})
provides a means of testing probabilistic models.
In this note the game-theoretic methodology
is extended to statistical models;
it will be demonstrated on the first-order scalar autoregressive model
\begin{equation}\label{eq:model}
  y_t
  =
  \alpha y_{t-1}
  +
  \epsilon_t,
  \quad
  t=1,2,\ldots,
\end{equation}
without the intercept term,
with constant $y_0$,
and with independent $N(0,1)$ innovations $\epsilon_t$.

We will be interested in procedures for computing,
for each $t=1,2,\ldots$,
a confidence interval $[l_t,u_t]$ for $\alpha$
given $y_0,\ldots,y_t$.
Let us fix a confidence level $1-\delta$,
and let $\alpha$ be the true parameter value.
The usual procedures are ``batch'',
in that they only guarantee
that $\alpha\in[l_t,u_t]$ with high probability
for a fixed $t$.
It is usually true that,
when they are applied sequentially,
the intersection $\cap_{t=1}^{\infty}[l_t,u_t]$ is empty with probability one.
Our goal is to guarantee that
\begin{equation}\label{eq:goal}
  \alpha\in\cap_{t=1}^{\infty}[l_t,u_t]
\end{equation}
with probability at least $1-\delta$.

Analogously to the usual classification of the limit theorems of probability theory
into ``strong'' (involving the conjunction over all $t$)
and ``weak'' (applicable to individual $t$),
let us call such confidence intervals \emph{strong}.
In particular,
confidence intervals satisfying (\ref{eq:goal})
with probability at least $1-\delta$
will be called \emph{strong $(1-\delta)$-confidence intervals}.
Accordingly, confidence intervals produced by the standard procedures
will be referred to as \emph{weak};
weak $(1-\delta)$-confidence intervals
satisfy $\alpha\in[l_t,u_t]$ with probability at least $1-\delta$
for each individual $t$.
(This probability is sometimes required to be precisely $1-\delta$,
but we will only consider the ``conservative'' definitions.)

To achieve the goal (\ref{eq:goal}),
for each possible value of the parameter $\alpha$
we construct random variables $S^{\alpha}_t$,
$t=0,1,\ldots$,
that form a nonnegative martingale
under the probability measure $\Prob_{\alpha}$
corresponding to the probabilistic model (\ref{eq:model})
with the given $\alpha$.
It will also be true that $S^{\alpha}_0=1$;
such sequences (nonnegative martingales starting from 1)
will be called \emph{martingale tests}.
We can then set
\begin{equation*}
  [l_t,u_t]
  =
  \left\{
    \alpha
    \st
    S^{\alpha}_t < 1/\delta
  \right\}
\end{equation*}
(assuming that the set on the right-hand side is an interval,
which it will be in our case).
The special case
\begin{equation*}
  \Prob_{\alpha}
  \left\{
    \sup_{t=0,1,\ldots}
    S^{\alpha}_t
    \ge
    1/\delta
  \right\}
  \le
  \delta
\end{equation*}
(due to Ville; see, e.g., \cite{ville:1939}, p.~100, or \cite{shafer/vovk:2001}, (2.12))
of Doob's inequality
shows that (\ref{eq:goal}) will indeed be true with probability
at least $1-\delta$.

\section{Derivation of strong confidence intervals}

If the true probability density of $y_t$
(conditional on the past) is
\begin{equation*}
  \frac{1}{\sqrt{2\pi}}
  \exp
  \left(
    -\frac
    {
      \left(
        y_t - \alpha\true y_{t-1}
      \right)^2
    }
    {
      2
    }
  \right)
\end{equation*}
and we want to reject the hypothesis
\begin{equation*}
  \frac{1}{\sqrt{2\pi}}
  \exp
  \left(
    -\frac
    {
      \left(
        y_t - \alpha y_{t-1}
      \right)^2
    }
    {
      2
    }
  \right),
\end{equation*}
the best, in many respects\footnote{%
  Cf., e.g., the nonnegativity of the Kullback--Leibler divergence,
  Neyman--Pearson lemma,
  and the optimality property of the probability ratio test in sequential analysis.},
martingale test
is the likelihood ratio sequence with the relative increments
\begin{multline*}
  \frac
  {
    \frac{1}{\sqrt{2\pi}}
    \exp
    \left(
      -\frac
      {
        \left(
          y_t - \alpha\true y_{t-1}
        \right)^2
      }
      {
        2
      }
    \right)
  }
  {
    \frac{1}{\sqrt{2\pi}}
    \exp
    \left(
      -\frac
      {
        \left(
          y_t - \alpha y_{t-1}
        \right)^2
      }
      {
        2
      }
    \right)
  }\\
  =
  \exp
  \left(
    \frac
    {
      (\alpha^2 - (\alpha\true)^2)
      y_{t-1}^2
      +
      2(\alpha\true-\alpha)
      y_{t-1} y_t
    }
    {
      2
    }
  \right).
\end{multline*}
The product over $t=1,\ldots,T$ is the martingale test itself:
\begin{equation}\label{eq:martingale-alpha-alpha-true}
  S_T^{\alpha,\alpha\true}
  =
  \exp
  \left(
    \frac
    {
      (\alpha^2 - (\alpha\true)^2)
      \Gamma_0
      +
      2(\alpha\true-\alpha)
      \Gamma_1
    }
    {
      2
    }
  \right),
\end{equation}
where
\begin{equation*}
  \Gamma_0
  =
  \sum_{t=1}^T
  y_{t-1}^2
\end{equation*}
and
\begin{equation*}
  \Gamma_1
  =
  \sum_{t=1}^T
  y_{t-1} y_t.
\end{equation*}
To get rid of the parameter $\alpha\true$,
let us integrate (\ref{eq:martingale-alpha-alpha-true})
over the probability distribution $N(\alpha,a^2)$ on the $\alpha\true$s:
\begin{multline*}
  S_T^{\alpha}
  =
  \frac{1}{\sqrt{2\pi}a}
  \int_{-\infty}^{\infty}
  \exp
  \left(
    \frac
    {
      (\alpha^2 - (\alpha\true)^2)
      \Gamma_0
      +
      2(\alpha\true-\alpha)
      \Gamma_1
    }
    {
      2
    }
  \right)\\
  \times
  \exp
  \left(
    -\frac{(\alpha-\alpha\true)^2}{2a^2}
  \right)
  \D\alpha\true\\
  =
  \frac{1}{\sqrt{2\pi}a}
  \int_{-\infty}^{\infty}
  \exp
  \left(
    -
    \left(
      \frac{1}{2a^2}
      +
      \frac{\Gamma_0}{2}
    \right)
    x^2
    +
    (\Gamma_1-\alpha\Gamma_0)
    x
  \right)
  \D x
\end{multline*}
(where I made the substitution $x=\alpha\true-\alpha$).
Now the formula
\begin{equation*}
  \int_{-\infty}^{\infty}
    \exp(-Ax^2+Bx)
  \D x
  =
  \sqrt{\frac{\pi}{A}}
  \exp
  \left(
    \frac{B^2}{4A}
  \right)
\end{equation*}
gives
\begin{equation}\label{eq:martingale-alpha}
  S_T^{\alpha}
  =
  \frac{1}{\sqrt{a^2\Gamma_0+1}}
  \exp
  \left(
    \frac{a^2}{2}
    \frac{(\Gamma_1-\alpha\Gamma_0)^2}{a^2\Gamma_0+1}
  \right).
\end{equation}

To find the confidence intervals corresponding to (\ref{eq:martingale-alpha}),
fix a confidence level $1-\delta$.
The $(1-\delta)$-confidence interval corresponding to (\ref{eq:martingale-alpha})
is defined as the set of $\alpha$s satisfying
\begin{equation*}
  \frac{1}{\sqrt{a^2\Gamma_0+1}}
  \exp
  \left(
    \frac{a^2}{2}
    \frac{(\Gamma_1-\alpha\Gamma_0)^2}{a^2\Gamma_0+1}
  \right)
  \le
  \frac{1}{\delta}.
\end{equation*}
Solving this in $\alpha$ gives the confidence interval
\begin{equation}\label{eq:confidence-interval-alpha}
  \left|
    \alpha - \frac{\Gamma_1}{\Gamma_0}
  \right|
  \le
  \sqrt
  {
    \frac{a^2\Gamma_0+1}{a^2\Gamma_0^2}
    \ln
    \frac{a^2\Gamma_0+1}{\delta^2}
  }.
\end{equation}

Notice that,
in the stationary case $\lvert\alpha\rvert<1$,
where $\Gamma_0$ has the order of magnitude $T$,
the size of the confidence interval (\ref{eq:confidence-interval-alpha})
is $O(\sqrt{\ln T / T})$ as $T\to\infty$.
This is worse that the usual iterated-logarithm behaviour
($O(\sqrt{\ln\ln T / T})$)
but agrees with \cite{nielsen:2005}, Theorem 2.5
(although the latter result is just an upper bound).
One can speculate that,
in the stationary case,
the $O(\sqrt{\ln\ln T / T})$ behaviour
will be recovered if the $N(\alpha,a^2)$ is replaced
by a probability distribution that is more concentrated around $\alpha$,
as in Ville's \cite{ville:1939} proof of the law of the iterated logarithm
(see also \cite{shafer/vovk:2001}, Chapter 5).

Most of the terms
in the confidence interval (\ref{eq:confidence-interval-alpha})
are familiar from the literature
(which, however, mainly covers the case of weak confidence intervals).
The centre $\frac{\Gamma_1}{\Gamma_0}$ of the interval
is just the least-squares estimate of $\alpha$ from the given sample.
The statistic
\begin{equation}\label{eq:tau}
  \tau_T
  =
  \left(
    \frac{\Gamma_1}{\Gamma_0} - \alpha
  \right)
  \sqrt{\Gamma_0}
  \approx
  \frac
  {
    \frac{\Gamma_1}{\Gamma_0} - \alpha
  }
  {
    \sqrt{\frac{a^2\Gamma_0+1}{a^2\Gamma_0^2}}
  }
\end{equation}
(for a fixed sample size $T$)
has been studied extensively.
In describing the known results I will follow \cite{chan/wei:1987}.
Mann and Wald \cite{mann/wald:1943} showed that $\tau_T$
is $N(0,1)$ asymptotically when $\lvert\alpha\rvert<1$.
Anderson \cite{anderson:1959} extended this to the case $\lvert\alpha\rvert>1$.
White \cite{white:1958} and Rao \cite{rao:1978} showed that,
in the case $\lvert\alpha\rvert=1$,
$\tau_T$ converges in distribution to
\begin{equation}\label{eq:unit-root}
  \frac{1}{2}
  \frac
  {
    W^2(1) - 1
  }
  {
    \sqrt
    {
      \int_0^1
        W^2(s)
      \D s
    }
  }
\end{equation}
where $W$ is a standard Brownian motion.

Suppose, for concreteness,
that (\ref{eq:tau}) is asymptotically $N(0,1)$.
The central asymptotic weak confidence interval for $\alpha$
based on the statistic given after the ``$\approx$'' in (\ref{eq:tau})
will be different from (\ref{eq:confidence-interval-alpha})
in that
\begin{equation}\label{eq:log-term}
  \sqrt
  {
    \ln
    \frac{a^2\Gamma_0+1}{\delta^2}
  }
  =
  \sqrt
  {
    2\ln
    \frac{1}{\delta}
    +
    \ln
    (a^2\Gamma_0+1)
  }
\end{equation}
will be replaced by the upper $\delta/2$-quantile of $N(0,1)$,
essentially by
\begin{equation*}
  \sqrt{2\ln\frac{2}{\delta}}
\end{equation*}
for a small $\delta$.
This is close to the first addend on the right-hand side of (\ref{eq:log-term}),
and so the second addend represents
the price that we are paying
for our confidence intervals being strong.

\section{Empirical results}

To test the test martingales (\ref{eq:martingale-alpha}) empirically,
I generated $y_0,\ldots,y_{1000}$ from the model (\ref{eq:model})
with $y_0=0$ and $\alpha=0.8,1$.
The case $\alpha=0.8$ illustrates the stationary behaviour ($\lvert\alpha\rvert<1$),
and the ``unit-root'' case $\alpha=1$ is intermediate
between the stationary and ``explosive'' ($\lvert\alpha\rvert>1$) behaviour.
Tables \ref{tab:coef0.8} and \ref{tab:coef1}
give the approximate weak central $99\%$-confidence intervals
based on the above approximations for $\tau_T$
(normal for $\alpha=0.8$ and (\ref{eq:unit-root}) for $\alpha=1$)
and the strong $99\%$-confidence intervals
computed from (\ref{eq:confidence-interval-alpha}).

\begin{table}[bt]
\centering
\begin{tabular}{c|c|c}
\hline
Type of the interval  & Confidence interval & Its width\\
\hline
Weak (approximate)    & $[0.736,0.837]$     & $0.101$\\
Strong                & $[0.716,0.857]$     & $0.141$\\
\hline
\end{tabular}
\caption{\label{tab:coef0.8}Weak and strong $99\%$-confidence intervals
  obtained for $T=1000$ and $\alpha=0.8$ (stationary case).
  The value of the constant $a$ is $0.1$.}
\end{table}

\begin{table}[bt]
\centering
\begin{tabular}{c|c|c}
\hline
Type of the interval  & Confidence interval & Its width\\
\hline
Weak  (approximate)   & $[0.982,1.003]$     & $0.022$\\
Strong                & $[0.977,1.010]$     & $0.033$\\
\hline
\end{tabular}
\caption{\label{tab:coef1}The analogue of Table \ref{tab:coef0.8}
  for $\alpha=1$ (unit root case).}
\end{table}

The intuition behind the value of $a$ in (\ref{eq:confidence-interval-alpha})
is that it should be of the same order of magnitude
as the expected width of the confidence interval
(since $a$ represents the order of magnitude
of the distance to the bulk of $\alpha\true$ that we are competing with).
It is taken as $0.1$ in the tables,
but the results will not be drastically different
if $a=1$, which is intuitively more ``neutral'',
is chosen: e.g.,
the width $0.141$ in Table \ref{tab:coef0.8}
would go up to $0.162$,
and the width $0.033$ in Table \ref{tab:coef1}
would go up to $0.037$.

\begin{figure}[bt]
  \centering
    \makebox{\includegraphics[width=\picturewidth]{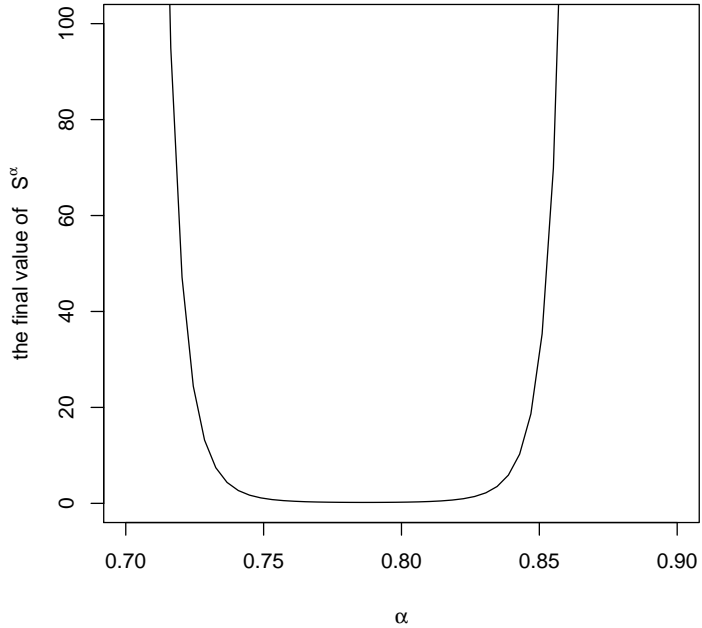}}
  \caption{\label{fig:coef0.8}The capital $S^{\alpha}_{1000}$
    achieved for various values of $\alpha$
    when the true coefficient $\alpha$ is $0.8$.}
\end{figure}

\begin{figure}[bt]
  \centering
    \makebox{\includegraphics[width=\picturewidth]{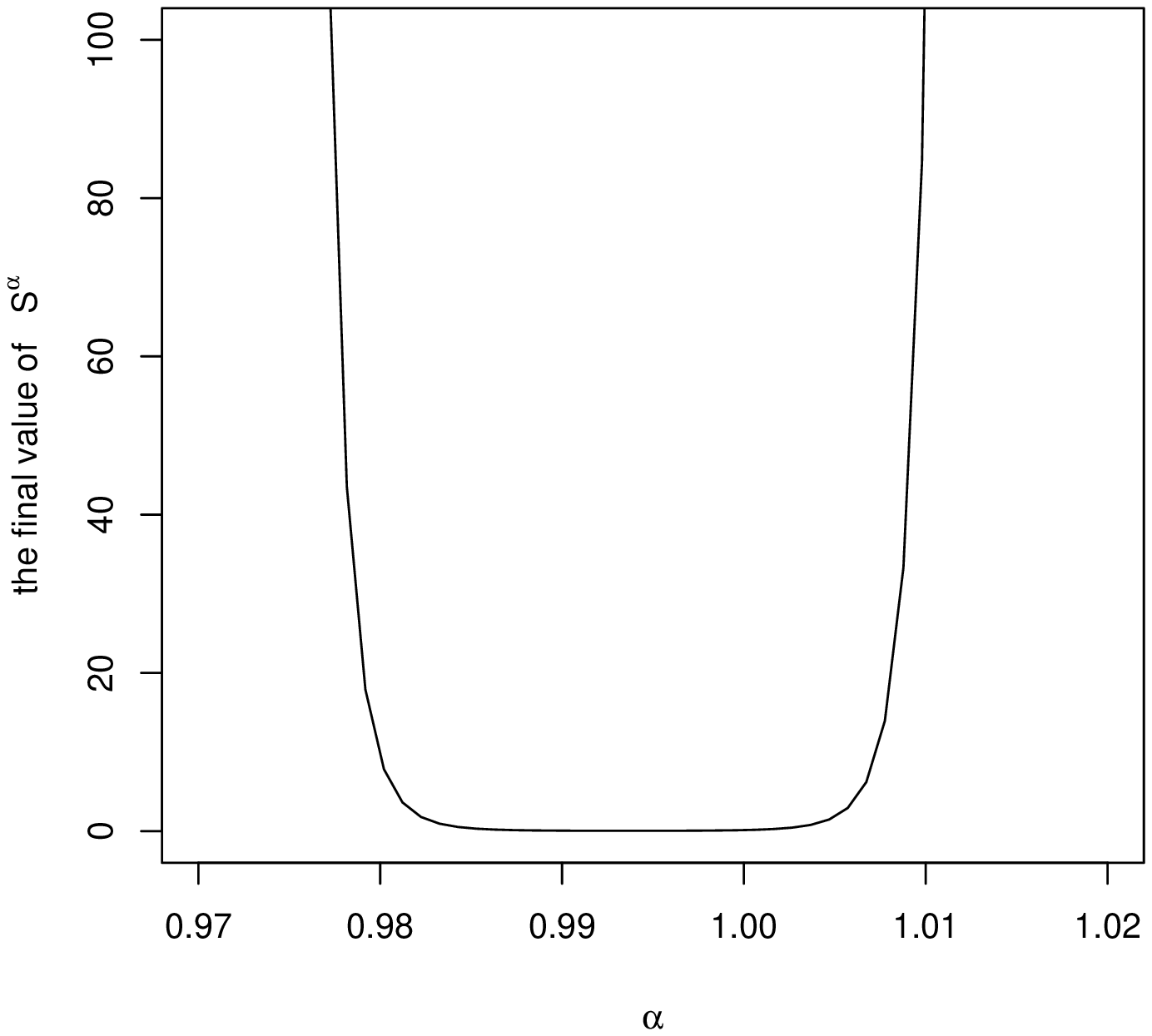}}
  \caption{\label{fig:coef1}The analogue of Figure \ref{fig:coef0.8}
    for $\alpha=1$.}
\end{figure}

Figures \ref{fig:coef0.8} and \ref{fig:coef1}
give the final values $S^{\alpha}_T$ for the same data set
and the same value of $a$,
$a=0.1$.

\section{Directions of further research}

These are some possible areas
in which the methods of martingale testing could be applied:
\begin{description}
\item[Online testing of statistical models.]
  When the strong confidence interval $[l_t,u_t]$ becomes empty,
  the statistical model can be rejected.
  Of course, efficient testing of statistical models
  will require different martingale tests:
  it will not be sufficient to consider, as in this note, different values of parameters
  as alternatives.
\item[Prediction.]
  In the simplest case,
  the prediction interval at step $t$
  might be computed as the union of the prediction intervals
  corresponding to all $\alpha\in[l_t,u_t]$.
\item[Alternative assumptions about innovations.]
  For example, the assumption that $\epsilon_t$ have zero medians
  (conditional on the past)
  might lead to feasible statistical procedures.
\end{description}

\subsection*{Acknowledgments}

I am grateful to Bent Nielsen, Clive Bowsher, David Hendry, and Jennifer Castle
for useful discussions.
This work was partially supported by EPSRC (grant EP/F002998/1),
MRC (grant G0301107),
and the Cyprus Research Promotion Foundation.


\begin{thebibliography}{1}
\bibitem{anderson:1959}
T\DOT{}~W\DOT{} Anderson.
\newblock On asymptotic distributions of estimates of parameters of stochastic
  difference equations.
\newblock {\em Annals of Mathematical Statistics}, 30:676--687, 1959.

\bibitem{chan/wei:1987}
N.~H. Chan and Ching~Zong Wei.
\newblock Asymptotic inference for nearly nonstationary {AR(1)} processes.
\newblock {\em Annals of Statistics}, 15:1050--1063, 1987.

\bibitem{mann/wald:1943}
H\DOT~B\DOT{} Mann and Abraham Wald.
\newblock On the statistical treatment of linear stochastic difference
  equations.
\newblock {\em Econometrica}, 11:173--220, 1943.

\bibitem{nielsen:2005}
Bent Nielsen.
\newblock Strong consistency results for least squares estimators in general
  vector autoregressions with deterministic terms.
\newblock {\em Econometric Theory}, 21:534--561, 2005.

\bibitem{rao:1978}
M\DOT~M\DOT{} Rao.
\newblock Asymptotic distribution of an estimator of the boundary parameter of
  an unstable process.
\newblock {\em Annals of Statistics}, 6:185--190, 1978.
\newblock Correction: 8:1403, 1980.

\bibitem{shafer/vovk:2001}
Glenn Shafer and \Vladimir{} Vovk.
\newblock {\em Probability and Finance: It's Only a Game!}
\newblock Wiley, New York, 2001.

\bibitem{ville:1939}
Jean Ville.
\newblock {\em Etude critique de la notion de collectif}.
\newblock Gauthier-Villars, Paris, 1939.

\bibitem{white:1958}
J\DOT~S\DOT{} White.
\newblock The limiting distribution of the serial correlation coefficient in
  the explosive case.
\newblock {\em Annals of Mathematical Statistics}, 29:1188--1197, 1958.
\end{thebibliography}
\end{document}